\documentclass[12pt,reqno,a4paper]{amsart}
\usepackage{graphicx}
\usepackage{amssymb,amscd,textcomp}
\usepackage{hyperref}

\begin{document}

\let\kappa=\varkappa
\let\eps=\varepsilon
\let\phi=\varphi

\def\Z{\mathbb Z}
\def\R{\mathbb R}
\def\C{\mathbb C}
\def\Q{\mathbb Q}

\def\OO{\mathcal O}
\def\CP{\C{\mathrm P}}
\def\RP{\R{\mathrm P}}
\def\conj{\overline}
\def\Beta{\mathrm{B}}
\def\p{\partial}

\renewcommand{\Im}{\mathop{\mathrm{Im}}\nolimits}
\renewcommand{\Re}{\mathop{\mathrm{Re}}\nolimits}
\newcommand{\codim}{\mathop{\mathrm{codim}}\nolimits}
\newcommand{\id}{\mathop{\mathrm{id}}\nolimits}
\newcommand{\Aut}{\mathop{\mathrm{Aut}}\nolimits}

\newtheorem*{thm}{Theorem}
\newtheorem{lem}{Lemma}

\theoremstyle{definition}
\newtheorem{rem}[lem]{Remark}
\newtheorem*{exm}{Example}

\title{Uniformization and Steinness}

\author[Nemirovski]{Stefan Nemirovski}
\address{%
Steklov Mathematical Institute, Moscow, Russia;\hfill\break
\strut\hspace{8 true pt} Fakult\"at f\"ur Mathematik, Ruhr-Universit\"at Bochum, Germany}
\email{stefan@mi.ras.ru}
\thanks{The first author was partly supported by SFB/TR 191 of the DFG and RFBR grant \textnumero 17-01-00592-a.
The second author was supported by a grant from the Natural Sciences and Engineering
Research Council of Canada.}
\author[Shafikov]{Rasul Shafikov}
\address{University of Western Ontario, London, Canada}
\email{shafikov@uwo.ca}
\begin{abstract}
It is shown that the unit ball in $\C^n$ is the only complex manifold
that can universally cover both Stein and non-Stein strictly pseudoconvex domains.
\end{abstract}


\maketitle

In this note we use methods from~\cite{NS2} to show that
the unit ball in $\C^n$ is the only simply connected complex manifold that can cover both 
Stein and non-Stein strictly pseudoconvex domains.

Here a strictly pseudoconvex domain is a relatively compact domain in a complex manifold 
such that its boundary admits a  $C^2$-smooth strictly plurisubharmonic defining 
function.

\begin{thm}
Let $Y$ be the universal cover of a Stein strictly pseudoconvex domain.
Suppose that $Y$ is not biholomorphic to the ball. Then any manifold covered by $Y$ does not
contain compact complex analytic subsets of positive dimension. In~particular, any other 
strictly pseudoconvex domain covered by $Y$ is Stein.
\end{thm}

Examples of strictly pseudoconvex domains covered by the ball in $\C^2$ which contain
compact complex curves (and hence are not Stein) may be found in~\cite{GKL}. It is 
well-known that the ball covers compact complex manifolds as well.

Recall also from~\cite{NS1,NS2} that a Stein strictly pseudoconvex domain 
is covered by the unit ball if and only if its boundary is everywhere locally 
CR-diffeomorphic to the unit sphere.

\smallskip 
The theorem will follow immediately from the two lemmas below.

\begin{lem}
Let $\pi:Y\to D$ be a covering of a complex manifold $D$ admitting a strictly plurisubharmonic
function $\phi:D\to\R$. If $A\subset Y$ is an analytic subset of positive dimension, 
then its projection $\pi(A)$ cannot lie in a compact subset in~$D$.
\end{lem}

\begin{rem}
The assumptions of the lemma are satisfied if $D$ is (an unramified domain over) a Stein manifold.
However, there exist examples of complex manifolds with strictly plurisubharmonic functions
but no non-constant holomorphic functions~\cite{F}. 
\end{rem}

\begin{proof}
Suppose that $\pi(A)$ is contained in a compact subset of $D$. Then there exists a sequence 
of points $x_n=\pi(y_n)$ such that $x_n\to x\in D$, $y_n\in A$, and 
$$
\sup\limits_{\pi(A)} \phi = \lim\limits_{n\to\infty}\phi(x_n)= \phi(x).
$$
Choose a convex coordinate neighbourhood $U\ni x$ and a strictly plurisubharmonic function
$\widetilde\phi$ on $\conj U$ such that $\widetilde\phi(x) = \phi(x)$ and 
$$
\widetilde\phi(\xi)< \phi(x) - \eps \quad \text{ for all } \xi\in \p U\cap \{\phi<\phi(x)\}
$$
for some $\eps>0$, see Fig.~\ref{accumulation}. 

\begin{figure}[htbp]
\begin{center}
\includegraphics[scale=0.6]{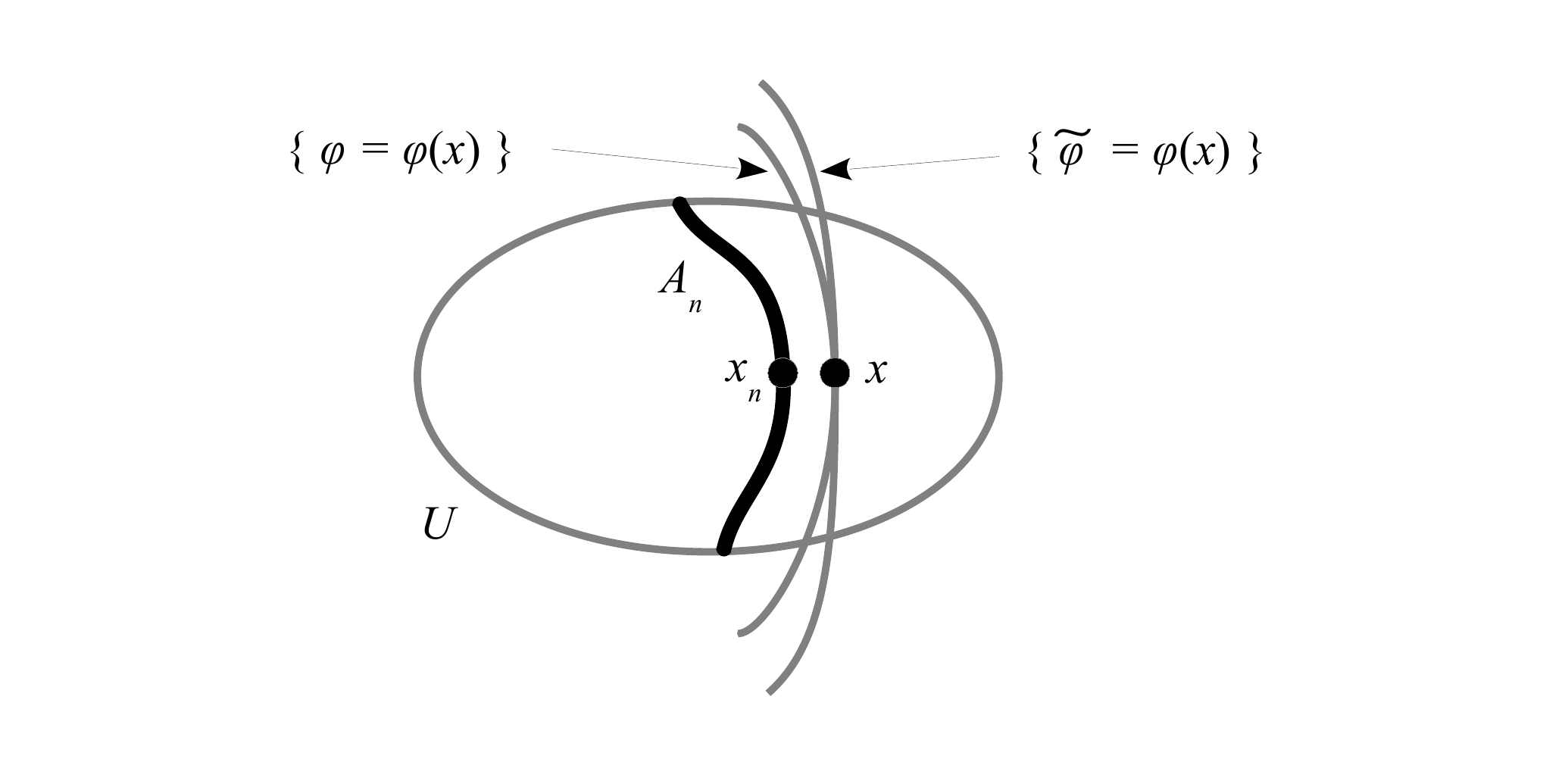}
\end{center}
\caption{Impossible analytic sets $A_n$.} 
\label{accumulation}
\end{figure}

For each $n\gg 1$, there is a local inverse $\psi_n$ to $\pi$ defined on $\conj U$ so that $\psi_n(x_n)=y_n$.
Set $A_n:= \psi_n^{-1}(A)$. This is a complex analytic subset of $U$ with boundary in $\p U$
and ${\phi|}_{A_n}\le \phi(x)$ by the choice of~$x$. Thus, ${\widetilde\phi|}_{\p A_n}<\phi(x) - \eps$
for all $n$ by construction, whereas $\widetilde\phi(x_n)\to\phi(x)$ as $n\to\infty$. This 
contradicts the maximum principle for plurisubharmonic functions on complex analytic
sets (see e.g.~\cite[\S 6.3]{C}), which proves the lemma.
\end{proof}

\begin{lem}
\label{lemma2}
Let $\pi:Y\to D$ be the universal covering of a strictly pseudoconvex domain by a complex manifold $Y$ 
that is not biholomorphic to the ball. Suppose that $\pi':Y\to M$  is a covering of a complex
manifold containing a connected compact complex analytic subset $B\Subset M$ of positive dimension.
Then $\pi\left(\pi'^{-1}(B)\right)$ is contained in a compact subset of~$D$. 
\end{lem}

\begin{rem}
In this lemma, $D$ does not need to be Stein.
\end{rem}

\begin{proof}
Let $\phi:\conj D\to (-\infty,0]$ be a plurisubharmonic defining function for $D$. 
Following~\cite[\S 2.3]{NS2}, consider the function $\psi$ on $M$ defined by
$$
\psi(x) := \left( \sup_{\pi'(y)=x} \phi\circ\pi(y) \right)^*,
$$
where $*$ denotes the upper semicontinuous regularisation. As shown in~\cite[\S 2.3]{NS2},
it follows from \cite[Corollary 2.3]{NS2} that $\psi$ is plurisubharmonic and {\it strictly\/} negative on~$M$. 
(It is explained in \cite[\S 3.2]{NS2} how to modify the proof of \cite[Corollary 2.3]{NS2} for non-Stein domains.)
By the maximum principle, 
$$
{\psi|}_B \equiv \mathrm{const.} <0.
$$
Hence, 
$$
\phi\circ\pi(y) \le \mathrm{const.} <0 \text{ for all } y\in\pi'^{-1}(B),
$$ 
and therefore $\pi\left(\pi'^{-1}(B)\right)$ is relatively compact in~$D$.
\end{proof}

\begin{rem}
The key point in the proof of Lemma~\ref{lemma2} is the application of \cite[Corollary 2.3]{NS2}.
That result is a consequence of \cite[Proposition 2.2]{NS2}, which is an extension of the 
well-known Wong--Rosay theorem~\cite{R,W} to universal coverings of strictly pseudoconvex
domains in complex manifolds.
\end{rem}

\end{document}